\documentclass[12pt]{article}

\usepackage{amsmath,amsfonts,amssymb,amsthm}
\usepackage{hyperref,ifthen,xifthen,iftex,rotating}
\hypersetup{
	pdfstartview=,
	pdfencoding=auto,
}

\usepackage[backend=biber,autopunct,sortcites,]{biblatex}

\hypersetup{linktocpage,unicode,pdfencoding=auto,pdfstartview=,}

\ifxetex
	\usepackage[no-math]{fontspec}
	\usepackage{unicode-math}
	\newcommand{\fakebold}{1}
	\newcommand{\scale}{1}
	\newcommand{\fakemathbold}{1}
	\newfontfamily{\bask}{ModernExtT-Theano.otf}
	\makeatletter
	\RenewDocumentCommand{\sum@}{}{\DOTSB\baskervillesum}
	\AtBeginDocument{%
		\RenewDocumentCommand{\sum}{}{\mathop{\sum@}\slimits@}%
	}
	\NewDocumentCommand{\baskervillesum}{}{%
		\mathchoice
		{\makebaskervillesum{1.8}}
		{\makebaskervillesum{1.2}}
		{\makebaskervillesum{1}}
		{\makebaskervillesum{0.7}}
	}
	\NewDocumentCommand{\makebaskervillesum}{m}{%
		\vcenter{\hbox{\scalebox{#1}{\bask Σ}}}%
	}
	\makeatletter
	\RenewDocumentCommand{\prod@}{}{\DOTSB\baskervilleprod}
	\AtBeginDocument{%
		\RenewDocumentCommand{\prod}{}{\mathop{\prod@}\slimits@}%
	}
	\NewDocumentCommand{\baskervilleprod}{}{%
		\mathchoice
		{\makebaskervilleprod{1.5}}
		{\makebaskervilleprod{1.2}}
		{\makebaskervilleprod{1}}
		{\makebaskervilleprod{0.7}}
	}
	\NewDocumentCommand{\makebaskervilleprod}{m}{%
		\vcenter{\hbox{\scalebox{#1}{\bask ∏}}}%
	}
	\newfontfamily\intfont{OldStandard-Italic.otf}
	\makeatletter
	\NewDocumentCommand{\int@}{}{\DOTSB\baskervilleint}
	\AtBeginDocument{%
		\RenewDocumentCommand{\int}{}{\mathop{\int@}\slimits@}%
	}
	\NewDocumentCommand{\baskervilleint}{}{%
		\mathchoice
		{\makebaskervilleint{1.5}}
		{\makebaskervilleint{1.2}}
		{\makebaskervilleint{1}}
		{\makebaskervilleint{0.7}}
	}
	\NewDocumentCommand{\makebaskervilleint}{m}{%
		\vcenter{\hbox{\scalebox{#1}{\intfont ∫}}}%
	}
	\defaultfontfeatures{Mapping=tex-text,Ligatures=Tex}
\else
	\usepackage[T1]{fontenc}
\fi

\DeclareMathOperator{\lcm}{lcm}

\newtheorem{theorem}{\textsc{\textbf{Theorem}}}
\theoremstyle{definition}
\newenvironment{definition}[1][]{\par\medskip\noindent \textsc{\textbf{{\ifthenelse{\isempty{#1}}{Definition.}{#1.}}}} \rmfamily}{\medskip}

\title{Memoir on Divisibility Sequences}
\author{\itshape Masum Billal}

\begin{document}
	\maketitle
		\begin{abstract}
			The purpose of this memoir is to discuss two very interesting properties of integer sequences. One is the law of apparition and the other is the law of repetition. Both have been extensively studied by mathematicians such as Ward, Lucas, Lehmer, Hall, etc. However, due to the lack of a proper survey in this area, many results have been rediscovered many decades later. This along with the necessity of the usefulness of such theory calls for a survey on this topic.
		\end{abstract}
	\section{Introduction}
	It is well known that we have $F_{m}\mid F_{n}$ for Fibonacci numbers $(F_{n})$ if $m\mid n$. In fact, we have $\gcd(F_{m},F_{n})=F_{\gcd(m,n)}$. \textcite{lucas_1878_1,lucas_1878_2,lucas_1878_3} and \textcite{lehmer30} generalized this property for Lucas sequence of the first kind $(U_{n})$ defined as
		\begin{align*}
			U_{n}
				& = \dfrac{\alpha^{n}-\beta^{n}}{\alpha-\beta}
		\end{align*}
	where $\alpha$ and $\beta$ are roots of $x^{2}-ax+b=0$ although under different conditions. They also establish the \textit{law of apparition} and the \textit{law of repetition}. The law of apparition is, if $\rho$ is the smallest index for which a prime $p$ divides $U_{\rho}$, then $p\mid U_{k}$ if and only if $\rho\mid k$. The law of repetition is, if $p^{\alpha}\|U_{\rho}$, then $p^{\alpha+\beta}\|U_{\rho p^{\beta}s}$ for $p\nmid s$.

	In this section, we discuss some basics. In \autoref{sec:elem}, we discuss properties of divisibility sequences in general. In \autoref{sec:lucas}, we will focus on the law of apparition for linear recurrences of order $k$. The reason we are so interested in the law of apparition becomes apparent once we have \autoref{thm:equiv}. In \autoref{sec:exp}, we investigate the law of repetition.
		\begin{definition}[Divisibility Sequence]
			An integer sequence $(a_{n})$ is a \textit{divisibility sequence} if $a_{m}\mid a_{n}$ whenever $m\mid n$. Some simple examples of divisibility sequences are $(n!),(\varphi(n)),(x^{n}-1),(F_{n})$.
		\end{definition}
	The term divisibility sequence was most likely used by \textcite{hall36} for the first time. Hall called a divisibility sequence $(a_{n})$ \textit{normal} if $a_{0}=0$ and $a_{1}=1$. We can actually assume that a divisibility sequence is normal without losing generality too much, as \textcite{hall36} has shown. In this memoir, we will be mostly concerned with the following stronger assumption.
		\begin{definition}[Strong Divisibility Sequence]
			An integer sequence $(a_{n})$ is a \textit{strong divisibility sequence} if $\gcd(a_{m},a_{n})=a_{\gcd(m,n)}$ for all positive integers $m$ and $n$. Some simple examples of strong divisibility sequences are $(x^{n}-1),(U_{n})$.
		\end{definition}
	 Although elliptic divisibility sequences are also divisibility sequences, we will not be focusing on that topic in this memoir. For elliptic divisibility sequences, the reader can consult \textcite{ward_1948}.
	 	\begin{definition}[Rank of Apparition]
	 		Let $m$ be a positive integer. If $\rho$ is the smallest index such that $m\mid a_{\rho}$, then $\rho$ is the \textit{rank of apparition} of $p$ in $(a_{n})$. For a prime $p$ and positive integer $e>1$, we denote the rank of apparition of $p^{e}$ by $\rho_{e}(p)$. If it is clear what the prime $p$ is, then we may only write $\rho_{e}$.
	 	\end{definition}

	 	\begin{definition}[Subsequence of Strong Divisibility Sequence]
	 		For a fixed positive integer $s$, the sequence $(c_{n})$ is a subsequence of $(a_{n})$ if
	 			\begin{align*}
	 				c_{n}
	 					& = \dfrac{a_{sn}}{a_{s}}
	 			\end{align*}
 			for all $n$.
	 	\end{definition}

 		\begin{definition}[Binomial Coefficients]
 			Let $n!_{a}$ denote the product of first $n$ terms of the strong divisibility sequence $(a_{n})$. Then the \textit{binomial coefficient} of $(a_{n})$ is
 				\begin{align*}
 					\binom{n}{k}_{a}
 						& = \dfrac{n!_{a}}{k!_{a}(n-k)!_{a}}
 				\end{align*}
 		\end{definition}
	\section{Elementary Properties}\label{sec:elem}
	We will first attempt to characterize strong divisibility sequences by its divisors. First, we see an analog of the law of repetition for strong divisibility sequences. A recent publication \textcite{billal_riasat_2021} discusses divisibility sequences and covers some of the results.
		\begin{theorem}
			Let $p$ be a prime and $\rho$ be the rank of apparition of $p$ in the strong divisibility sequence $(a_{n})$. Then $p\mid a_{k}$ if and only if $\rho\mid k$.
		\end{theorem}

		\begin{theorem}
			Let $m$ be a positive integer and the prime factorization of $m$ be
				\begin{align*}
					m
						& = \prod_{i=1}^{r}p_{i}^{e_{i}}
				\end{align*}
			If the rank of apparition of $p_{i}^{e_{i}}$ in $(a_{n})$ is $\rho_{e_{i}}(p_{i})$, then the rank of apparition of $m$ is
				\begin{align*}
					\rho
						& = \lcm(\rho_{e_{1}}(p_{1}),\ldots,\rho_{e_{r}}(p_{r}))
				\end{align*}
		\end{theorem}
	We have the first necessary and sufficient condition for a divisibility sequence $(a_{n})$ to be a strong divisibility sequence due to \textcite{ward36}.
		\begin{theorem}\label{thm:equiv}
			Let $(a_{n})$ be a divisibility sequence. Then $(a_{n})$ is a strong divisibility sequence is equivalent to the condition that for a prime $p$ and positive integer $e$, $p^{e}\mid a_{k}$ if and only if $\rho_{e}(p)\mid k$.
		\end{theorem}
	\textcite{ward55_2} proves the following result. \textcite{nowicki} essentially rediscovers the same result.
		\begin{theorem}
			Let $(a_{n})$ be an integer sequence. Then $(a_{n})$ is a strong divisibility sequence if and only if there exists an integer sequence $(b_{n})$ such that
				\begin{align*}
					a_{n}
						& = \prod_{d\mid n}b_{d}
				\end{align*}
			where $\gcd(b_{m},b_{n})=1$ whenever $m\nmid n$ and $n\nmid m$.
		\end{theorem}

		\begin{definition}[LCM Sequence]
			This new sequence $(b_{n})$ associated with $(a_{n})$ is the \textit{lcm sequence} of $(a_{n})$. It can be thought of as a generalization of cyclotomic polynomials $\Phi_{n}(x)$ of $x^{n}-1$.
		\end{definition}

		\begin{theorem}
			Let $(a_{n})$ be a strong divisibility sequence and $(b_{n})$ is the lcm sequence of $(a_{n})$. Then
				\begin{align*}
					\lcm(a_{1},\ldots,a_{n})
						& = b_{1}\cdots b_{n}
				\end{align*}
		\end{theorem}

		\begin{theorem}
			The lcm sequence $(b_{n})$ of a strong divisibility sequence $(a_{n})$ is given by
				\begin{align*}
					b_{n}
						& = \dfrac{\lcm(a_{1},\ldots,a_{n})}{\lcm(a_{1},\ldots,a_{n-1})}\\
						& = \dfrac{a_{n}\prod_{\substack{p_{i},p_{j}\mid n\\i\neq j}}a_{\frac{n}{p_{i}p_{j}}}}{\prod_{p_{i}\mid n}a_{n/p_{i}}\prod_{\substack{p_{i},p_{j},p_{k}\mid n\\i\neq j\neq k}}a_{\frac{n}{p_{i}p_{j}p_{j}}}}\\
						& = \dfrac{a_{n}}{\lcm(a_{n/p_{1}},\ldots,a_{n/p_{r}})}
				\end{align*}
			where $p_{1},\ldots,p_{r}$ are distinct prime factors of $n$.
		\end{theorem}

		\begin{theorem}
			Let $(a_{n})$ be an integer sequence. Then $(a_{n})$ is a strong divisibility sequence if and only if for a positive integer $m>1$ and positive integers $k,l$, we have $m\mid a_{k},m\mid a_{l}$ if and only if $m\mid a_{\gcd(k,l)}$.
		\end{theorem}
	A corollary is the following.
		\begin{theorem}\label{thm:onerank}
			A divisibility sequence $(a_{n})$ is a strong divisibility sequence if and only any positive integer $m>1$ assumes only one rank of apparition.
		\end{theorem}

		\begin{theorem}
			If an integer sequence $(u_{n})$ has the property that $\gcd(u_{pn},u_{qn})=u_{n}$ for distinct primes $p,q$ and positive integers $n$, let us say that $(u_{n})$ has property P. Then both the strong divisibility sequence $(a_{n})$ and its lcm sequence $(b_{n})$ have the property P.
		\end{theorem}

		\begin{theorem}
			If $(a_{n})$ is a divisibility sequence and $\gcd(a_{pn},a_{qn})=a_{n}$ for distinct primes $p$ and $q$, then $\gcd(a_{m},a_{n})=1$ if $\gcd(m,n)=1$.
		\end{theorem}

		\begin{theorem}
			A necessary and sufficient condition that an integer sequence $(a_{n})$ is a strong divisibility sequence is that
				\begin{align*}
					\gcd(a_{pn},a_{qn})
						& = a_{n}
				\end{align*}
			for all distinct primes $p,q$ and positive integers $n$.
		\end{theorem}
	We have the analogous of Legendre's theorem for strong divisibility sequences.
		\begin{theorem}
			Let $(a_{n})$ be a strong divisibility sequence and $p$ be a prime. Then
				\begin{align*}
					\nu_{p}(n!_{a})
						& = \sum_{i\geq 1}\left\lfloor{\dfrac{n}{\rho_{i}(p)}}\right\rfloor
				\end{align*}
		\end{theorem}

		\begin{theorem}
			The binomial coefficients of a strong divisibility sequence are integers.
		\end{theorem}
	\section{Lucasian Sequences}\label{sec:lucas}
	In this section, we will see the connection between linear recurrent and divisibility sequences. Some of the results will make use of abstract algebra when it seems convenient to do so. But we will mostly concern ourselves with integer sequences since analogous results usually extend to the appropriate field.
		\begin{definition}[Linear Recurrent Sequence]
			A \textit{linear recurrent sequence} of order $k$ is defined as
				\begin{align}
					u_{n+k}
						& = c_{k-1}u_{n+k-1}+\ldots+c_{0}u_{n}\label{eqn:linrec}
				\end{align}
		\end{definition}
	We are interested in $(u_{n})$ when the coefficients $c_{0},\ldots,c_{k-1}$ are integers. We can easily extend the definition over a field $\mathbb{F}$. The polynomial associated with $(u_{n})$ in \autoref{eqn:linrec} is the \textit{characteristic polynomial} of $u$ which is
		\begin{align*}
			f(x)
				& = x^{k}-c_{k-1}x^{k-1}-\ldots-c_{0}
		\end{align*}
	Denote the discriminant of $f$ by $\mathfrak{D}(f)$. If it is clear what $f$ is, we may write $\mathfrak{D}$ only.
		\begin{definition}[Lucasian Sequence]
			An integer sequence $(u_{n})$ is \textit{Lucasian} if $u$ is both a linear recurrent sequence and a divisibility sequence. \textcite{ward37,ward55_2} called such sequences ``Lucasian" \textit{in honor of the french mathematician E. Lucas who first systematically studied a special class of such sequences}.
		\end{definition}

		\begin{definition}[Null Divisor]
			A positive integer $n$ is a \textit{null divisor} of the Lucasian sequence $(u_{n})$ if $n\mid u_{m}$ for all $m\geq n_{0}$. If $(u_{n})$ has no null divisor other than $1$, then $(u_{n})$ is \textit{primary}. $d$ is a \textit{proper null divisor} of $(u_{n})$ if $d$ divides neither the initial terms $u_{0},\ldots,u_{k-1}$ nor the coefficients $c_{0},\ldots,c_{k-1}$. If $d$ is not a proper null divisor, then it is a \textit{trivial null divisor}.
		\end{definition}

		\begin{definition}[Generator]
			Define the polynomial $f_{i}$ as $f_{0}(x)=0$ and
				\begin{align*}
					f_{r}
						& = x^{r}-c_{r-1}x^{r-1}-\ldots-c_{0}
				\end{align*}
			Then the polynomial
				\begin{align*}
					\mathfrak{u}(x)
						& = u_{0}f_{k-1}(x)+\ldots+u_{k-1}f_{0}(x)
				\end{align*}
			is called the \textit{generator} of $(u_{n})$. If
				\begin{align*}
					\Delta(\mathfrak{u})
						& =
							\begin{vmatrix}
								u_{0} & \ldots & u_{k-1}\\
								u_{1} & \ldots & u_{k}\\
								\vdots & \ddots & \vdots\\
								u_{k-1} & \ldots & u_{2k-2}
							\end{vmatrix}
				\end{align*}
			then we have
				\begin{align*}
					\Delta(\mathfrak{u})
						& = (-1)^{k(k-1)/2}\mathfrak{R}(u(x),f(x))
				\end{align*}
			where $\mathfrak{R}(f(x),g(x))$ is the \textit{resultant} of two polynomials $f$ and $g$.
		\end{definition}

		\begin{definition}[Index]
			Let $\nu_{n}(a)$ be the largest non-negative integer $k$ such that $n^{k}\mid a$ but $n^{k+1}\nmid a$. If $G$ is the largest null divisor of $(u_{n})$, then for a proper null prime divisor $p$, $\nu_{p}(G)$ is the \textit{index} of $p$ in $(u_{n})$.
		\end{definition}

		\begin{definition}[Period and Numeric]
			Consider the Lucasian sequence $(u_{n})$ modulo $m$. Let $\rho$ be the least positive index such that
				\begin{align*}
					U_{\rho}
						& \equiv 0\pmod{m}\\
						& \vdots\\
					U_{\rho+k-2}
					& \equiv 0\pmod{m}\\
					U_{\rho+k-1}
						& \equiv 1\pmod{m}
				\end{align*}
			Then $\rho$ is a \textit{period} of $(u_{n})$ modulo $m$ because
				\begin{align*}
					u_{n+\rho}
						& \equiv u_{n}\pmod{m}
				\end{align*}
			for all $n\geq n_{0}$. The number of non-periodic terms of $(u_{n})$ modulo $m$ is the \textit{numeric}. We say that $(u_{n})$ is \textit{periodic} modulo $m$ and $(u_{n})$ is \textit{purely periodic} modulo $m$ if the numeric $n_{0}=0$. On the other hand, $\tau$ is a \textit{restricted period} of $(u_{n})$ modulo $m$ if $\tau$ is the least positive integer for which
				\begin{align*}
					U_{\tau}
						& \equiv 0\pmod{m}\\
						& \vdots\\
					U_{\tau+k-2}
						& \equiv 0\pmod{m}
				\end{align*}
			In this case, $u_{n+\tau}\equiv Au_{n}\pmod{m}$ for some $m\nmid A$ and all $n\geq n_{0}'$. This $A$ is called the \textit{multiplier} of $(u_{n})$ modulo $m$. The value of this multiplier $A$ depends on $\tau$.
		\end{definition}

		\begin{definition}[R-sequence]
			Let $(u_{n})$ be a Lucasian sequence with an irreducible polynomial $f$. If $\alpha_{1},\ldots,\alpha_{k}$ are the roots of $f$, then
				\begin{align*}
					U_{n}(f)
						& = \prod_{i<j}\left(\dfrac{\alpha_{i}^{n}-\alpha_{j}^{n}}{\alpha_{i}-\alpha_{j}}\right)
				\end{align*}
			is the \textit{R-sequence} associated with $(u_{n})$. We simply write $U_{n}$ if it is clear what $f$ is. Then $(U_{n})$ is a Lucasian sequence. The case $k=2$ gives us the classical Lucas sequence of the first kind. R-sequences are of particular importance because Lucasian sequences seem to be either R-sequences themselves or divisors of R-sequences. Moreover, the consideration of R-sequence gives us further insight into the determination of the law of apparition.
		\end{definition}

		\begin{definition}[Period of Polynomial]
			Let $f$ be a polynomial irreducible modulo $p$. Then the smallest positive integer $n$ for which
				\begin{align*}
					x^{n}
						& \equiv 1\pmod{p,f(x)}
				\end{align*}
			is the \textit{period of $f$ modulo }$p$. For two polynomials $h(x)$ and $g(x)$, we write
				\begin{align*}
					g(x)
						& \equiv h(x)\pmod{m,f(x)}
				\end{align*}
			if
				\begin{align*}
					g(x)-h(x)
						& = f(x)q(x)+m\cdot r(x)
				\end{align*}
			for some polynomials $q$ and $r$.
		\end{definition}
	\textcite{hall36} states the following easily derived results.
		\begin{theorem}\label{thm:redchar}
			Let $(u_{n})$ be a normal Lucasian sequence with characteristic polynomial $f$ such that the prime $p$ does not divide the discriminant $\mathfrak{D}(f)$. If
				\begin{align*}
					f(x)
						& \equiv f_{1}(x)\cdots f_{s}(x)\pmod{p}
				\end{align*}
			is the factorization of $f$ modulo $p$ into irreducible polynomials $f_{1},\ldots,f_{s}$ of degree $k_{1},\ldots,k_{s}$ and $\rho$ is the least period of $(u_{n})$ modulo $p$, then
				\begin{align*}
					\rho
						& \mid \lcm(p^{k_{1}}-1,\ldots,p^{k_{s}}-1)
				\end{align*}
		\end{theorem}
	Due to \autoref{thm:redchar}, we can turn our attention primarily to the case when $f$ is irreducible modulo the prime $p$.
		\begin{theorem}
			Let $(u_{n})$ be a normal Lucasian sequence. If $\rho$ is a rank of apparition and $\tau$ is a restricted period of $(u_{n})$ modulo the prime $p$ respectively, then $\rho\mid\tau$.
		\end{theorem}

		\begin{theorem}
			Let $(u_{n})$ be a normal Lucasian sequence and $\tau$ be its restricted period modulo the prime $p$. If $p\mid n$, then $\tau\mid n$.
		\end{theorem}
	Note that this result is slightly stronger than the typical result that the rank of apparition $\rho\mid n$ if $p\mid u_{n}$ since $\rho\mid\tau$ but the converse is not always true. \textcite{ward38} proves the following generalized result.
		\begin{theorem}
			Let $\mathfrak{O}$ be a commutative ring and $(u_{n})$ be a Lucasian sequence with elements in $\mathfrak{O}$. Moreover, $\mathfrak{A}$ is an ideal of $\mathfrak{O}$ such that no divisor of $\mathfrak{A}$ is a null divisor of $(u_{n})$. Then if $(u_{n})$ is periodic modulo $\mathfrak{A}$, the minimal restricted period of $(u_{n})$ modulo $\mathfrak{A}$ exists and divides every other restricted period of $(u_{n})$. This minimal restricted period divides the period of $(u_{n})$ modulo $\mathfrak{A}$. Furthermore, the multipliers of $(u_{n})$ modulo $\mathfrak{A}$ are relatively prime to $\mathfrak{A}$ and forms a group with respect to multiplication modulo $\mathfrak{A}$.
		\end{theorem}

		\begin{theorem}
			Let $\mathfrak{O}$ be a ring and $(u_{n})$ be a sequence of $\mathfrak{O}$ and $\mathfrak{A}$ be an ideal such that $(u_{n})$ is periodic modulo $\mathfrak{A}$ but no divisor of $\mathfrak{A}$ is a null divisor of $(u_{n})$. If $\rho$ is the least period and $\tau$ is the restricted period of $(u_{n})$ modulo $\mathfrak{A}$, then the multipliers of $(u_{n})$ form a cyclic group of order $\rho/\tau$. Furthermore, the multiplier dependent on $\tau$ is a  of this group.
		\end{theorem}
	The concept of the rank of apparition is almost the same as the rank of apparition of strong divisibility sequences for Lucasian sequences. However, unlike strong divisibility sequences, it is possible that sometimes $(u_{n})$ may have more than one rank of apparition modulo $\mathfrak{A}$. For this reason, we can probably redefine the rank of apparition of $\mathfrak{A}$ in the following way. We call $\rho$ a rank of apparition of $\mathfrak{A}$ in $(u_{n})$ for the ring $\mathfrak{O}$ if
		\begin{align*}
			u_{\rho}
				& \equiv 0\pmod{\mathfrak{A}}\\
			\iff u_{d}
				& \not\equiv 0\pmod{\mathfrak{A}}
		\end{align*}
	for any divisor $d$ of $\rho$. With this connection, one of our primary interests is knowing when the set of the rank of apparitions is finite. Note that, when we consider such a set of ranks of apparition, we can actually consider a rank of apparition $\delta$ a duplicate of the rank of apparition $\rho$ if $\rho\mid\delta$. The obvious reason being that the ranks covered by $\delta$ are already covered by $\rho$. In this regard, we have the following result.
		\begin{theorem}
			Let $\mathfrak{A}$ be a divisor of the Lucasian sequence $(u_{n})$ such that $(u_{n})$ is periodic modulo $\mathfrak{A}$. Then a necessary and sufficient condition that $\mathfrak{A}$ has a finite set of ranks of apparition in $(u_{n})$ is that all the ranks divide the restricted period of $(u_{n})$ modulo $\mathfrak{A}$.
		\end{theorem}

		\begin{theorem}
			Let $(u_{n})$ be a Lucasian sequence and $\mathfrak{A}$ be a divisor of $(u_{n})$ such that $(u_{n})$ is purely periodic modulo $\mathfrak{A}$. Then $\mathfrak{A}$ only has a finite set of ranks and each rank divides the restricted period of $(u_{n})$ modulo $\mathfrak{A}$.
		\end{theorem}
	Let $m$ be a positive integer that does not divide the coefficient $c_{0}$ of $u$ and $\mathfrak{S}_{m}$ denote the set of all ranks of apparition of $(u_{n})$ modulo $m$. We readily have the following result.
		\begin{theorem}\label{thm:finrank}
			The set $\mathfrak{S}_{m}$ consists of all multiples of a finite set of rank of apparition $\rho_{1},\ldots,\rho_{s}$ such that
			\begin{align*}
				u_{\rho_{i}}
				& \equiv 0\pmod{m}\\
				\iff u_{d}
				& \not\equiv 0\pmod{m}
			\end{align*}
			for any $d\mid\rho_{i}$ and $\rho_{i}\nmid\rho_{j}$.
		\end{theorem}
	The finite set in \autoref{thm:finrank} is called the \textit{ranks of apparition} of $(u_{n})$ modulo $m$. We can actually consider $(u_{n})$ modulo $m$ using a \textit{single unified rank of apparition }$\rho$ where $\rho=\lcm(\rho_{1},\ldots,\rho_{s})$. The places of apparition of $m$ in $(u_{n})$ are periodic modulo $\rho$ and $\rho\mid\tau$ where $\tau$ is the restricted period of $(u_{n})$.
		\begin{theorem}
			Let $(u_{n})$ be a normal Lucasian sequence of order $k$ and $\mathfrak{l}=\lcm(1,\ldots,k)$. Then $p^{k}(p^{\mathfrak{l}}-1)$ is a period of $(u_{n})$ modulo $p$.
		\end{theorem}

		\begin{theorem}
			Let $(u_{n})$ be a Lucasian sequence of order $k$ with characteristic polynomial $f(x)$ and $p$ be a prime. If $p\mid u_{p}$, then $p\mid\mathfrak{D}(f)$ or $p\mid c_{0}$.
		\end{theorem}

		\begin{theorem}
			Let $p$ be a null divisor of a normal Lucasian sequence $(u_{n})$, then $p$ divides both $\Delta(\mathfrak{u})$ and $\mathfrak{D}(f)$ where $\mathfrak{u}$ is the generator and $f(x)$ is the characteristic polynomial of $u$ respectively.
		\end{theorem}

		\begin{theorem}
			A sufficient condition that the Lucasian sequence $(u_{n})$ is primary is that $\gcd(\Delta(\mathfrak{u}),\mathfrak{D}(f))=1$ where $\mathfrak{u}$ is the generator and $f$ is the characteristic polynomial of $(u_{n})$ respectively.
		\end{theorem}

		\begin{theorem}
			Let $p$ be a null prime divisor of a Lucasian sequence $(u)$ such that the coefficients are relatively prime. If $\mathfrak{u}$ is the generator of $(u_{n})$, then $\nu_{p}(\Delta(\mathfrak{A}))$ is the index of $p$ in $(u_{n})$.
		\end{theorem}

		\begin{theorem}
			A subsequence of a normal Lucasian sequence can have no prime null divisor that is not a possible null divisor of $(u_{n})$ itself.
		\end{theorem}

		\begin{theorem}
			Let $(u_{n})$ be a primary Lucasian sequence of order $k$ such that the characteristic polynomial has no repeated roots, the coefficients are relatively prime and $\mathfrak{l}=\lcm(1,\ldots,k)$. Then
				\begin{align*}
					u_{p}^{\mathfrak{l}}
						& \equiv 1\pmod{p}
				\end{align*}
			for large enough $p$.
		\end{theorem}

		\begin{theorem}
			Let $(u_{n})$ be a Lucasian sequence with characteristic polynomial $f$, $(U_{n})$ be the associated R-sequence and $p$ be a prime such that $p\nmid \mathfrak{D}(f)$. Then every rank of apparition of $p$ in $(U_{n})$ is a rank of apparition in $(u_{n})$.
		\end{theorem}
	Next, we have a generalization of the law of apparition given by Lucas.
		\begin{theorem}\label{thm:genrank}
			Let $(u_{n})$ be a Lucasian sequence of order $k$ with characteristic polynomial $f$ irreducible modulo $p$ and $\lambda$ be the period of $f$ modulo $p$. If $k$ has the prime factorization
				\begin{align*}
					k
						& = q_{1}^{e_{1}}\cdots q_{s}^{e_{s}}
				\end{align*}
			then the ranks of apparition of $p$ in $(U_{n})$ are divisors of the elements of a subset of
				\begin{align*}
					\{\rho(k/q_{1}),\ldots,\rho(k/q_{s})\}
				\end{align*}
			where $\rho(s)=\lambda/\gcd(\lambda,p^{s}-1)$. Thus, $p$ has at most $k$ distinct ranks of apparition and the single unified rank of $p$ divides
				\begin{align*}
					\rho\left(\dfrac{k}{q_{1}\cdots q_{s}}\right)
				\end{align*}
		\end{theorem}
	A corollary is the following.
		\begin{theorem}\label{thm:lucdivex}
			Any Lucasian sequence with an irreducible characteristic polynomial of order $k$ where $k$ is a prime power has only one rank of apparition and hence, is a strong divisibility sequence.
		\end{theorem}

		\begin{theorem}
			The Lucasian sequence $(u_{n})$ is not a strong divisibility sequence if it has an irreducible characteristic polynomial and the ranks of apparitions are in the set
				\begin{align*}
					\{\rho(k/q_{1}),\ldots,\rho(k/q_{r})\}
				\end{align*}
			for $1<r<s$ where $q_{1},\ldots,q_{s}$ are the distinct prime divisors of $k$.
		\end{theorem}

		\begin{theorem}
			The prime $p$ is a null divisor of the Lucasian sequence $(U_{n})$ if and only if $p$ divides the last two coefficients $c_{1}$ and $c_{0}$ of the characteristic polynomial $f$ of $(u_{n})$.
		\end{theorem}
	\section{The Law of Repetition}\label{sec:exp}
	We say that an integer sequence $(a_{n})$ has the \textit{law of repetition} if for any positive integer $n$ and a prime divisor $p$ of $a_{n}$ such that $p\nmid s$,
		\begin{align*}
			\nu_{p}(a_{nk})
				& = \nu_{p}(a_{n})+\nu_{p}(k)
		\end{align*}
	holds.
		\begin{theorem}\label{thm:expdiv}
			Let $(a_{n})$ be an integer sequence with the law of repetition. Then $(a_{n})$ is also a strong divisibility sequence.
		\end{theorem}

		\begin{proof}
			For positive integers $m$ and $n$, let $g=\gcd(m, n),m=gu,n=gv$ where $\gcd(u,v)=1$ and $h=\gcd(a_{m},a_{n})$. We will show that $h=a_{g}$. First, consider that $p$ is a prime divisor of $g$. If $p^{e}\|a_{g}$,
				\begin{align*}
					\nu_{p}(h)
						& = \min\left(\nu_{p}(a_{gu}),\nu_{p}(a_{gv})\right)\\
						& = \nu_{p}(a_{g})+\min(\nu_{p}(u),\nu_{p}(v))
				\end{align*}
			Since $\gcd(u,v)=1$, $p$ cannot divide both $u$ and $v$. Therefore, either $\nu_{p}(u)$ or $\nu_{p}(v)$ is $0$ and $\min(\nu_{p}(u),\nu_{p}(v))=0$. This gives us $\nu_{p}(h)=\nu_{p}(a_{g})$ for all prime divisor $p$ of $g$. Next, assume that $p$ is a prime divisor of $h$ and $p^{e}\|h$. Then $p^{e}\mid a_{m}$ and $p^{e}\mid a_{n}$. More specifically, $p^{e}\|a_{gu}$ or $p^{e}\|a_{gv}$ must hold. Again, by definition $\nu_{p}(a_{gu})=\nu_{p}(a_{g})+\nu_{p}(u)$ and $\nu_{p}(a_{gv})=\nu_{p}(a_{g})+\nu_{p}(v)$. Since both $p\mid u$ and $p\mid v$ cannot hold, so $p^{e}\|a_{gu}$ or $p^{e}\|a_{gv}$ must hold. Then $p^{e}\|a_{g}$ holds for all $p^{e}\|h$. Thus, we must have $h=a_{g}$.
		\end{proof}
	By \autoref{thm:expdiv}, any sequence with the law of repetition has a corresponding lcm sequence $(b_{n})$. The next result characterizes when a strong divisibility sequence has the law of repetition.
		\begin{theorem}\label{thm:expchar}
			Let $(a_{n})$ be a strong divisibility sequence, $(b_{n})$ be the lcm sequence of $(a_{n})$ and $\rho$ be the rank of apparition of prime $p$ in $(a_{n})$. Then $(a_{n})$ has the law of repetition if and only if for any positive integers $n$ and $m>1$ such that $p\nmid m$, $p\|b_{\rho p^n}$ but $p\nmid b_{\rho p^nm}$.
		\end{theorem}

		\begin{proof}
			First, we will prove the if part. Since $(a_{n})$ is a strong divisibility sequence, $p\mid a_{k}$ if and only if $\rho\mid k$. By assumption, $(a_{n})$ has law of repetition. If $p^\alpha\|a_\rho$, then $p^{\alpha+1}\|a_{\rho p}$.
				\begin{align*}
					a_{\rho p}
						& = \prod_{d\mid \rho p}b_d\\
					\nu_{p}(a_{\rho p})
						& = \nu_{p}\left(\prod_{d\mid\rho p}b_d\right)
				\end{align*}
			If $d<\rho$, then $p\nmid a_d$ so $p\nmid b_d$. Thus,
				\begin{align*}
					\nu_{p}(a_{\rho p})
						& = \nu_{p}\left(\prod_{d\mid p}b_{\rho d}\right)\\
						& = \nu_{p}(b_\rho)+\nu_{p}(b_{\rho p})\\
					\alpha+1
						& = \alpha+\nu_{p}(b_{\rho p})
				\end{align*}
			So, $\nu_{p}(a_{\rho p})=1$ and $p\mid b_{\rho p}$. By induction, we can see that $p$ not only divides $b_{\rho p^i}$ for $i\in\mathbb{N}$, more precisely, $p\|b_{\rho p^i}$. Next, assume that $p^{\alpha+u}\|a_{n}$ for some positive integer $n=\rho p^{u}m$ where $p\nmid m$. From the law of repetition and the argument above,
				\begin{align*}
					\nu_{p}(a_{n})
						& = \nu_{p}(a_{\rho p^{u}m})\\
						& = \nu_{p}(a_{\rho})+\nu_{p}\left(\prod_{d\mid p^um}b_{\rho d}\right)\\
						& = \alpha+\nu_{p}\left(\prod_{d\mid p^u}b_{\rho d}\right)+\nu_{p}\left(\prod_{\substack{d\mid p^u\\e\mid m\\e>1}}b_{\rho de}\right)
				\end{align*}
			Since $\nu_{p}(a_{n})=\nu_{p}(a_{\rho p^{u}m})=\nu_{p}(a_{\rho})+u$,
				\begin{align*}
					\alpha+u
						& = \alpha+\sum_{i=1}^u\nu_{p}(b_{\rho p^i})+\nu_{p}\left(\prod_{i=1}^u\prod_{\substack{e\mid m\\e>1}}b_{\rho p^ie}\right)\\
						& = \alpha+u+\nu_{p}\left(\prod_{i=1}^u\prod_{\substack{e\mid m\\e>1}}b_{\rho p^ie}\right)\\
						& = \alpha+u+\sum_{i=1}^u\sum_{\substack{e\mid m\\e>1}}\nu_{p}(b_{\rho p^ie})
				\end{align*}
			From this, we have that $\nu_{p}(b_{\rho p^ie})=0$ for $1\leq i\leq u$ and $e\mid m$ if $e>1$. In other words, $p\mid b_k$ if and only if $k=\rho p^u$ for some non-negative integer $u$.

			For the only if part, we have that $(a_{n})$ is a strong divisibility sequence such that $p\| b_{\rho p^u}$ but $p\nmid b_{\rho p^um}$ for $m>1$. Let $n$ be a positive integer such that $n=\rho p^um$ and $p^\alpha\|a_\rho$.
			\begin{align*}
				\nu_{p}(a_{n})
					& = \nu_{p}(a_{\rho p^um})\\
					& = \nu_{p}\left(\prod_{d\mid \rho p^um}b_d\right)\\
					& = \nu_{p}(a_\rho)+\nu_{p}\left(\prod_{d\mid p^um}b_{\rho d}\right)\\
			\end{align*}
		Now, separate the sum into two parts based on whether the index has a divisor of $m$ greater than $1$.
			\begin{align*}
				\nu_{p}(a_{n})
					& = \nu_{p}(a_{\rho})+\sum_{d\mid p^u}\nu_{p}(b_{\rho d})+\sum_{d\mid p^u}\sum_{\substack{e\mid m\\e>1}}\nu_{p}(b_{\rho de})\\
					& = \alpha+\sum_{i=1}^u\nu_{p}(b_{\rho p^i})+0\\
					& = \alpha+\sum_{i=1}^u1\\
					& = \alpha+u
			\end{align*}
			This proves the theorem.
		\end{proof}
	A corollary of \autoref{thm:expchar} is the following.
		\begin{theorem}\label{thm:lte}
			Let $(a_{n})$ be a sequence with the law of repetition and $(b_{n})$ be the lcm sequence of $(a_{n})$. If $m$ and $n$ are distinct positive integers, then $\gcd(b_{m},b_{n})>1$ if and only if $m/n$ is a prime power. More precisely, $p$ is a prime divisor of $\gcd(b_{m},b_{n})$ if and only if $m/n=p^{s}$ for some non-negative integer $s$.
		\end{theorem}
	\printbibliography
\end{document}